\documentclass[11pt, a4paper]{amsart}
\usepackage[latin1]{inputenc}
\usepackage{amssymb,amsmath,amsthm}
\usepackage{enumerate}
\usepackage{nicefrac}

\usepackage{cases}
\usepackage{hyperref}

\usepackage[mathcal]{eucal}

\newcommand{\Z}{{\mathbb Z}}

\def\coker{\operatorname{Coker}}

\newtheorem{Thm}{Theorem}[section]
\newtheorem{Prop}[Thm]{Proposition}
\newtheorem{Cor}[Thm]{Corollary}
\newtheorem{Lem}[Thm]{Lemma}
\newtheorem*{PCThm}{P\'{o}lya--Carlson Theorem}

\theoremstyle{definition}

\newtheorem{Ex}[Thm]{Example}

\newtheorem*{Defi*}{Definition}

\theoremstyle{remark}

\newtheorem{Rmk}[Thm]{Remark}

\numberwithin{equation}{section}

\title[Coincidence Reidemeister zeta functions]{P\'olya--Carlson
  dichotomy for coincidence Reidemeister zeta functions via profinite
  completions}



\author[A.\ Fel'shtyn]{Alexander Fel'shtyn}

\address{Instytut Matematyki, Uniwersytet Szczecinski, Szczecin, Poland}
\email{alexander.felshtyn@usz.edu.pl}

\author[B.\ Klopsch]{Benjamin Klopsch$^*$}
\address{Heinrich-Heine-Universit\"{a}t D\"{u}sseldorf, Mathematisches
  Institut, Universit\"{a}tsstr. 1, 40225 D\"{u}sseldorf, Germany}
\email{klopsch@math.uni-duesseldorf.de}

\thanks{$^*$ Corresponding author \\ The research was funded by the
  Deutsche Forschungsgemeinschaft (DFG, German Research Foundation) --
  380258175; and it was funded by the Narodowy Centrum Nauki of Poland
  (NCN) (grant No.2016/23/G/ST1/04280 (Beethoven 2)).  The research
  was also partially conducted in the framework of the DFG-funded
  research training group ``GRK 2240: Algebro-Geometric Methods in
  Algebra, Arithmetic and Topology''.}

\begin{document}
	
\begin{abstract}
  We consider coincidence Reidemeister zeta functions for tame
  endomorphism pairs of nilpotent groups of finite rank, shedding new
  light on the subject by means of profinite completion techniques.

  In particular, we provide a closed formula for coincidence
  Reidemeister numbers for iterations of endomorphism pairs of
  torsion-free nilpotent groups of finite rank, based on a weak
  commutativity condition, which derives from simultaneous
  triangularisability on abelian sections.  Furthermore, we present
  results in support of a P\'olya--Carlson dichotomy between
  rationality and a natural boundary for the analytic behaviour of the
  zeta functions in question.
\end{abstract}
	
\maketitle

%%%%%
	
\section{Introduction}  \label{sec:introduction}

In classical topological fixed point theory, Reidemeister numbers
arise as homotopy invariants associated to iterates of a continuous
self-map of a connected compact polyhedron.  Passing to the
fundamental group of the polyhedron, the Reidemeister numbers admit an
algebraic treatment in terms of twisted conjugacy classes; compare
\cite{Fe00}.  In this paper we take a group-theoretic point of view,
inspired by, but otherwise largely independent of the topological
origins of the subject.

Let $G$ be a group and let $\varphi, \psi \colon G \to G$ be
endomorphisms of~$G$.  Elements $x,y \in G$ are said to be
\emph{$(\varphi,\psi)$-twisted conjugate} to one another if there
exists $g \in G$ such that $x = (g \varphi)^{-1} y (g \psi)$.  We
observe that this sets up an equivalence relation on~$G$; the
corresponding equivalence classes are called
\emph{$(\varphi,\psi)$-twisted conjugacy classes} or
\emph{$(\varphi,\psi)$-coincidence Reidemeister classes}.  We denote
by $\mathcal{R}(\varphi,\psi)$ the set of all
$(\varphi,\psi)$-coincidence Reidemeister classes of~$G$, and
$R(\varphi,\psi) = \lvert \mathcal{R}(\varphi,\psi) \rvert$ is called
the \emph{$(\varphi,\psi)$-coincidence Reidemeister number}.

We call the pair $(\varphi,\psi)$ of endomorphisms \emph{tame} if the
Reidemeister numbers $R(\varphi^n,\psi^n)$ are finite for all $n \in
\mathbb{N}$.  For such a tame pair of endomorphisms we define the    
\emph{$(\varphi,\psi)$-coincidence Reidemeister zeta function}
\[
Z_{\varphi,\psi}(s) = \exp \left( \sum_{n=1}^\infty
  \frac{R(\varphi^n,\psi^n)}{n} s^n \right),
\]
where $s$ denotes a complex variable.  This zeta function can be
regarded as an analogue of the Hasse--Weil zeta function of an
algebraic variety over a finite field or the Artin--Mazur zeta
function of a continuous self-map of a topological space.  In the
theory of dynamical systems, the coincidence Reidemeister zeta
function counts the synchronisation points of two maps, i.e.\ the
points whose orbits intersect under simultaneous iteration of two
endomorphisms; see \cite{Mi13}, for instance.  For
$\psi = \mathrm{id}_G$, we recover ordinary $\varphi$-twisted
conjugacy classes (or $\varphi$-Reidemeister classes) and the ordinary
$\varphi$-Reidemeister zeta function that were studied, for instance
in~\cite{Fe91,Fe00,Wo01,DeDu15,FeLe15,DeTeBu18,FeZi20,FeTrZi20}.  For
short, we set
$\mathcal{R}(\varphi) = \mathcal{R}(\varphi,\mathrm{id}_G)$ and
$R(\varphi) = R(\varphi,\mathrm{id}_G)$; if $(\varphi,\mathrm{id}_G)$
is tame, we say that $\varphi$ is tame and write
$Z_\varphi(s) = Z_{\varphi,\mathrm{id}_G}(s)$.  From a group-theoretic
point of view, the conincidence Reidemeister zeta function is a
natural generalisation of the ordinary Reidemeister zeta function.

\begin{Rmk}
  (1) The anti-isomorphism $G \to G$, $x \mapsto x^{-1}$ sets up a
  one-to-one correspondence between $(\varphi,\psi)$-twisted and
  $(\psi,\varphi)$-twisted conjugacy classes.  Consequently,
  $(\varphi,\psi)$ is tame if and only if $(\psi,\varphi)$ is tame,
  and $Z_{\varphi,\psi}(s) = Z_{\psi,\varphi}(s)$ for tame pairs.

  (2) If $\psi \colon G \to G$ is an automorphism, it is easy to check
  that $\mathcal{R}(\varphi,\psi) = \mathcal{R}(\psi^{-1}\varphi)$ and
  consequently $R(\varphi,\psi) = R(\psi^{-1}\varphi)$.  Moreover, if
  $(\varphi,\psi)$ is tame and if $\varphi \psi = \psi \varphi$ (in
  some situations a weak concept of commutativity suffices, as we will
  see), then $\psi^{-1}\varphi$ is tame and
  $Z_{\varphi,\psi}(s) = Z_{\psi^{-1} \varphi}(s)$.
\end{Rmk}

In order to illustrate these concepts by means of a simple example,
namely the infinite cyclic group, and for later use we record the
following easy fact.

\begin{Lem} \label{lem:img-psi-central}
  Let $\varphi, \psi \colon G \to G$ be endomorphisms of a group $G$
  such that the image of $\psi$ is central in~$G$, i.e.\
  $G\psi \subseteq \mathrm{Z}(G)$.  Then
  $H_{\varphi,\psi} = \{ (g \varphi)^{-1} (g \psi) \mid g \in G \}$ is
  a subgroup of~$G$ and $\mathcal{R}(\varphi,\psi) = H_{\varphi,\psi}
  \backslash G$; consequently, $R(\varphi,\psi) = \lvert G :
  H_{\varphi,\psi} \rvert$.
\end{Lem}

\begin{Ex}\label{ex1}
  Let $G = \mathbb{Z}$ be the infinite cyclic group, written
  additively, and let
  \[
  \varphi \colon G \to G, \quad x \mapsto d_\varphi x \qquad \text{and} \qquad
  \psi \colon G \to G, \quad x \mapsto d_\psi x
  \]
  for $d_\varphi, d_\psi \in \mathbb{Z}$.  According to the lemma, for every $n
  \in \mathbb{N}$, we have 
  \[
  R(\varphi^{\, n},\psi^{\, n}) =
  \begin{cases}
    \lvert d_\psi^{\, n} - d_\varphi^{\, n} \rvert & \text{if $d_\varphi^{\, n}
      \not = d_\psi^{\, n}$,} \\
    \infty & \text{otherwise}.
  \end{cases}
  \] 
  Consequently, $(\varphi,\psi)$ is tame precisely when
  $\lvert d_\varphi \rvert \not = \lvert d_\psi \rvert$ and, in this case, 
  \[
  Z_{\varphi,\psi}(s) = \frac{1 - d_2 s}{1 - d_1 s} \qquad \text{where
    $d_1 = \max \{ \lvert d_\varphi \rvert, \lvert d_\psi \rvert\}$
    and
    $d_2 = \frac{d_\varphi d_\psi}{d_1}$.}
  \]
\end{Ex}

This simple example or at least special cases of it are known.  The
aim of the current paper is to generalise the example -- as far as
possible -- to torsion-free nilpotent groups of finite (Pr\"ufer)
rank; the notion of rank and other relevant concepts are recalled in
Section~\ref{subsec:nilp}.  Our approach via profinite completions
offers new techniques and sheds fresh light on the subject.  We prove
a number of results, some of which provide generalisations from
ordinary Reidemeister classes to coincidence Reidemeister classes and,
at the same time, from finitely generated nilpotent groups to
nilpotent groups of finite rank.

In the following theorem we summarise some of our results.  We write
$\mathbb{P}$ for the set of all rational primes; for
$p \in \mathbb{P}$, the field of $p$-adic numbers is denoted by
$\mathbb{Q}_p$, the ring of $p$-adic integers by $\mathbb{Z}_p$, and
the $p$-adic absolute value (as well as its unique extension to the
algebraic closure $\overline{\mathbb{Q}}_p$) by
$\lvert \cdot \rvert_p$.  The absolute value on $\mathbb{C}$ is
denoted by $\lvert \cdot \rvert_\infty$.

\begin{Thm} \label{thm:main-result-nilpotent} Let
  $\varphi, \psi \colon G \to G$ be a tame pair of endomorphisms of a
  torsion-free nilpotent group~$G$ of finite Pr\"ufer rank.  Let $c$
  denote the nilpotency class of~$G$ and, for $1 \le k \le c$, let
  $\varphi_k ,\psi_k \colon G_k \to G_k$ denote the induced
  endomorphisms of the torsion-free abelian factor groups
  $G_k=\overline{\gamma}_k(G) / \overline{\gamma}_{k+1}(G)$ of finite
  rank, $d_k\geq 1$ say, that arise from the isolated lower central
  series~\eqref{equ:isolated-lcs} of~$G$.  Then the following hold.

  \smallskip
  
  \textup{(1)} For each $n \in \mathbb{N}$, there is a bijection
  between the set $\mathcal{R}(\varphi^n, \psi^n)$ of
  $(\varphi,\psi)$-coincidence Reidemeister classes and the cartesian
  product $\prod_{k=1}^c \mathcal{R}(\varphi_k^{\, n}, \psi_k^{\, n})$;
  consequently,
  \[
    R(\varphi^n, \psi^n) = \prod_{k=1}^c R(\varphi_k^{\, n},
    \psi_k^{\, n}) \qquad \text{for $n \in \mathbb{N}$.}
  \]
  
  \textup{(2)} For $1 \le k \le c$, let
  \[
    \varphi_{k,\mathbb{Q}}, \psi_{k,\mathbb{Q}} \colon
    G_{k,\mathbb{Q}} \to
    G_{k,\mathbb{Q}}
  \]
  denote the extensions of $\varphi_k, \psi_k$ to the divisible hull
  $G_{k,\mathbb{Q}} = \mathbb{Q} \otimes_\mathbb{Z} G_k \cong
  \mathbb{Q}^{d_k}$ of~$G_k$.  Suppose that each pair of endomorphisms
  $\varphi_{k,\mathbb{Q}}, \psi_{k,\mathbb{Q}}$ is simultaneously
  triangularisable.  Let $\xi_{k,1}, \ldots, \xi_{k,d_k}$ and
  $\eta_{k,1}, \ldots, \eta_{k,d_k}$ be the eigenvalues of
  $\varphi_{k,\mathbb{Q}}$ and $\psi_{k,\mathbb{Q}}$ in a fixed
  algebraic closure of the field~$\mathbb{Q}$, including
  multiplicities, ordered so that, for $n \in \mathbb{N}$, the
  eigenvalues of
  $\varphi_{k,\mathbb{Q}}^{\, n} - \psi_{k,\mathbb{Q}}^{\, n}$ are
  $\xi_{k,1}^{\, n} - \eta_{k,1}^{\, n}, \ldots, \xi_{k,d_k}^{\, n} -
  \eta_{k,d_k}^{\,n}$.  Set
  $L_k = \mathbb{Q}(\xi_{k,1}, \ldots, \xi_{k, d_k}, \eta_{k,1},
  \ldots, \eta_{k, d_k})$; for each $p \in \mathbb{P}$ fix an
  embedding $L_k \hookrightarrow \overline{\mathbb{Q}}_p$ and choose
  an embedding $L_k \hookrightarrow \mathbb{C}$.

  Then there exist subsets $I_k(p) \subseteq \{1,\ldots,d_k\}$, for
  $p \in \mathbb{P}$, such that the following hold.
  \begin{enumerate}[\rm (i)]
  \item For each $p \in \mathbb{P}$, the polynomials
    $\prod_{i \in I_k(p)} (X - \xi_{k,i})$ and
    $\prod_{i \in I_k(p)} (X - \eta_{k,i})$ have coefficients in
    $\mathbb{Z}_p$; in particular,
    $\lvert \xi_{k,i} \rvert_p , \lvert \eta_{k,i} \rvert_p \le 1$ for
    $i \in I_k(p)$.
  \item For each $n \in \mathbb{N}$,
  \begin{equation} \label{equ:key-formula-2}
    R(\varphi_k^{\, n},\psi_k^{\, n}) = \prod_{p \in \mathbb{P}}
    \prod_{i \in I_k(p)} \lvert \xi_{k,i}^{\, n} - \eta_{k,i}^{\, n}
    \rvert_p^{\, -1} = \prod_{i=1}^{d_k }\lvert \xi_{k,i}^{\, n} -
    \eta_{k,i}^{\, n} \rvert_\infty \cdot \prod_{p \in \mathbb{P}}
    \prod_{i \not \in I_k(p)} \lvert \xi_{k,i}^{\, n} - \eta_{k,i}^{\,
      n} \rvert_p;
    \end{equation}
    as this number is a positive integer,
    $\lvert \xi_{k,i}^{\, n} - \eta_{k,i}^{\, n} \rvert_p = 1$ for
    $1 \le i \le d_k$ for almost all~$p \in \mathbb{P}$.
  \end{enumerate}

  \smallskip
  
  \textup{(3)} Suppose that, for each $k \in \{1, \ldots, c\}$, the
  primes $p$ that contribute non-trivial factors
  $\prod_{i \not\in I_k(p)} \lvert \xi_{k,i}^{\, n} - \eta_{k,i}^{\,
    n} \rvert_p \neq 1$ to the product on the far right-hand side of
  \eqref{equ:key-formula-2} form a finite subset
  $\mathbf{P}_k \subseteq \mathbb{P}$ and that
  $\lvert \xi_{k,i} \rvert_\infty \neq \lvert \eta_{k,i}
  \rvert_\infty$ for $1 \le i \le d_k$.
  
  Then the coincidence Reidemeister zeta function
  $Z_{\varphi,\psi}(s)$ is either a rational function or it has a
  natural boundary at its radius of convergence.  Furthermore, the
  latter occurs if and only if
  $\lvert \xi_{k,i} \rvert_p= \lvert \eta_{k,i} \rvert_p$ for some
  $k \in \{1,\ldots,c\}$, $p\in\mathbf{P}_k$ and $i \not\in I_k(p)$.
\end{Thm}

For a full discussion and several intermediate results, some of which
are also of independent interest, we refer to the main text.  The
paper is organised as follows.

Section~\ref{subsec:nilp} puts results of Roman'kov~\cite{Ro11} on
ordinary Reidemeister classes of finitely generated torsion-free
nilpotent groups in the more general perspective of coincidence
Reidemeister classes of torsion-free nilpotent groups of finite rank.
Of particular interest is Proposition~\ref{Nil}.  In
Section~\ref{subsec:prof-compl} we establish a natural bijection
between the set of coincidence Reidemeister classes of a pair
$(\varphi,\psi)$ of endomorphisms of an almost abelian group $G$ and
the corresponding set for the induced pair
$(\widehat{\varphi},\widehat{\psi})$ of endomorphisms of the profinite
completion~$\widehat{G}$; see Proposition~\ref{profinite-bijection}.
As explained and illustrated there, this closes little gaps in the
literature.

In Section~\ref{sec:explicit-formula} we derive a closed formula for
the sequence of coincidence Reidemeister numbers for iterations of
tame endomorphism pairs of abelian groups of finite rank.  Our
approach is via profinite completions and sheds new light on a similar
formula of Miles~\cite{Mi08}, which was established in a rather
different way, namely based upon techniques from commutative algebra.
Corollary~\ref{cor:reduction} provides a reduction to torsion-free
groups, which are then duly dealt with in Proposition~\ref{formula}.

In Section~\ref{sec:polya-carlson} we present, in analogy to works of
Bell, Miles, Ward~\cite{BeMiWa14} and Byszewski, Cornelissen~\cite[\S
5]{ByCo18}, results in support of a P\'olya--Carlson dichotomy between
rationality and a natural boundary for the analytic behaviour of the
coincidence Reidemeister zeta functions of torsion-free nilpotent
groups of finite rank.  The underlying ideas are already present in
the abelian case which is covered in Theorem~\ref{abel}.  In contrast
to previous results, the P\'olya-Carlson dichotomy established in
Theorem~\ref{thm:main-result-nilpotent} applies for the first time to
a large class of non-abelian, not necessarily finitely generated
groups.  The simple, but important Example~\ref{exa:non-comm-endos}
illustrates a new phenomenon that occurs for coincidence Reidemeister
zeta functions associated to non-commuting pairs of endomorphisms.
The paper concludes with the proof of
Theorem~\ref{thm:main-result-nilpotent} already stated above.

%%%%%

\section{Nilpotent groups and profinite completions}

\subsection{Nilpotent groups} \label{subsec:nilp} A group $G$ has
\emph{finite (Pr\"ufer) rank} if there exists an integer $r$ such that
every finitely generated subgroup of $G$ can be generated by $r$
elements; the least such $r$ is called the \emph{(Pr\"ufer) rank}
of~$G$.  It is known that $G$ is a torsion-free nilpotent group of
finite rank if and only if there is some integer $n$ such that $G$ is
isomorphic to a subgroup of $\mathrm{Tr}_1(n,\mathbb{Q})$, the group
of all upper uni-triangular matrices over~$\mathbb{Q}$; compare
\cite[Cor.\ to Thm.~2.5]{We73}.

Let $G$ be a torsion-free nilpotent group of class $c = c(G)$.  As
usual, let
$G = \gamma_1(G) \ge \gamma_2(G) \ge \ldots \ge \gamma_c(G) \ge
\gamma_{c+1}(G) = 1$ denote the lower central series.  The elements of
finite order in a nilpotent group form a subgroup, the \emph{torsion
  subgroup} of ~$G$.  For $1 \le i \le c+1$, define
$\overline{\gamma}_i(G)$ to be the subgroup of $G$ such that
$\overline{\gamma}_i(G)/\gamma_i(G)$ is the torsion subgroup of
$G/\gamma_i(G)$.  Then it is routine to check that
\begin{equation} \label{equ:isolated-lcs} G = \overline{\gamma}_1(G)
  \ge \overline{\gamma}_2(G) \ge \ldots \ge \overline{\gamma}_c(G) \ge
  \overline{\gamma}_{c+1}(G) = 1
\end{equation}
is a descending central series of fully invariant subgroups of~$G$
such that each factor group
$\overline{\gamma}_i(G) / \overline{\gamma}_{i+1}(G)$ is torsion-free,
for $1 \le i \le c$; for instance, see \cite[Chap.~11,
Lem.~1.8]{Pa77}.
% For completeness, we sketch the argument.  It is straightforward that
% each $\overline{\gamma}_i(G)$ is fully invariant in~$G$ and that
% $G = \overline{\gamma}_1(G) \supseteq \ldots \supseteq
% \overline{\gamma}_{c+1}(G) = 1$.  Clearly, each factor group
% $\overline{\gamma}_i(G) / \overline{\gamma}_{i+1}(G)$ is torsion-free.
% It remains to show that
% $\overline{\gamma}_i(G) / \overline{\gamma}_{i+1}(G)$ is abelian and,
% in fact, central in $G / \overline{\gamma}_{i+1}(G)$.  In both cases
% it is enough to check the assertion for finitely generated subgroups
% of~$G$; hence without loss of generality we take $G$ to be finitely
% generated.  By construction,
% $\overline{\gamma}_i(G) / \overline{\gamma}_{i+1}(G)$ is
% central-by-finite; thus Schur's theorem~\cite[\S 10.1.4]{Ro96} implies
% that $\overline{\gamma}_i(G) / \overline{\gamma}_{i+1}(G)$ is
% finite-by-abelian, hence abelian.  Now
% $\gamma_i(G) \overline{\gamma}_{i+1}(G) / \overline{\gamma}_{i+1}(G)$
% is central in $G/ \overline{\gamma}_{i+1}(G)$ and of finite index in
% the finitely generated torsion-free abelian group
% $\overline{\gamma}_i(G) / \overline{\gamma}_{i+1}(G)$.  This implies
% that $\overline{\gamma}_i(G) / \overline{\gamma}_{i+1}(G)$ is central
% in $G / \overline{\gamma}_{i+1}(G)$.
% This establishes all claimed properties of the
% series~\eqref{equ:isolated-lcs}; ...
We refer to the series~\eqref{equ:isolated-lcs} as the \emph{isolated
  lower central series} of~$G$.  If, in addition, $G$ has finite rank
then each factor group
$\overline{\gamma}_i(G) / \overline{\gamma}_{i+1}(G)$ is a
torsion-free abelian group of rank $d_i$, say, and hence isomorphic to
a subgroup of the additive vector group $\mathbb{Q}^{d_i}$ that
contains a $\mathbb{Q}$-basis.

The following proposition puts results of Roman'kov~\cite{Ro11} for
ordinary twisted conjugacy classes in finitely generated nilpotent
groups in a more general perspective.
 
\begin{Prop} \label{Nil}
  Let $G$ be a torsion-free nilpotent group, of class~$c$ and with
  isolated lower central series~\eqref{equ:isolated-lcs}.  Let
  $\varphi, \psi \colon G \to G$ be endomorphisms and, for
  $1 \le i \le c$, let
  $\varphi_i ,\psi_i \colon \overline{\gamma}_i(G) /
  \overline{\gamma}_{i+1}(G) \to \overline{\gamma}_i(G) /
  \overline{\gamma}_{i+1}(G)$
  denote the induced endomorphisms of the torsion-free abelian factor
  groups arising from~\eqref{equ:isolated-lcs}.

  \textup{(1)} Suppose that $G$ has finite rank and that
  $R(\varphi,\psi) < \infty$.  Then, for $1 \le i \le c$, the
  coincidence set
  \[
  \mathrm{Coin}(\varphi_i,\psi_i) = \big\{ x. \overline{\gamma}_{i+1}(G)
  \in \overline{\gamma}_i(G) / \overline{\gamma}_{i+1}(G) \mid
  \big( x. \overline{\gamma}_{i+1}(G) \big) \varphi_i =
  \big( x. \overline{\gamma}_{i+1}(G) \big) \psi_i \big\}
  \]
  is a singleton, i.e.\ equal to $\{ 1. \overline{\gamma}_{i+1}(G) \}$.

  \textup{(2)} Suppose that, for $1 \le i \le c$, the coincidence set
  $\mathrm{Coin}(\varphi_i,\psi_i)$ is a singleton.  Then there is a
  bijection between $\mathcal{R}(\varphi, \psi)$ and
  $\prod_{i=1}^c \mathcal{R}(\varphi_i, \psi_i)$.  In particular,
  $R(\varphi, \psi) < \infty$ if and only if
  $R(\varphi_i, \psi_i) < \infty$ for $1 \le i \le c$, and in this
  case $R(\varphi, \psi) = \prod_{i=1}^c R(\varphi_i, \psi_i)$.
\end{Prop}

\begin{proof}
  It suffices to check that the arguments developed in~\cite{Ro11} can
  be adapted so that they work in the more general situation
  considered here.  By induction, we may assume that $c \ge 1$ and
  that the analogous assertions hold for the endomorphisms of
  $\overline{\varphi}, \overline{\psi} \colon G /
  \overline{\gamma}_c(G) \to G / \overline{\gamma}_c(G)$
  that are induced by $\varphi$ and~$\psi$.
  
  Clearly, two elements $x,y \in G$ can be $(\varphi,\psi)$-twisted
  conjugate only if their images in $G / \overline{\gamma}_c(G)$ are
  $(\overline{\varphi},\overline{\psi})$-twisted conjugate to one
  another.  Furthermore, we observe that, for any central subgroup $A$
  of $G$,
  \[
  L_{G,\varphi,\psi}(A) = \{ a \in A \mid \exists g \in G : (g
  \varphi) a = g \psi \} = \{ a \in A\mid \exists g \in G : a = (g
  \varphi)^{-1} (g \psi) \}
  \]
  is a subgroup of~$A$.  If, in addition, $A$ is $\varphi$- and
  $\psi$-invariant, then
  \[
    A / L_{G,\varphi,\psi}(A) \to \{ [a] \in \mathcal{R}(\varphi,\psi)
    \mid a \in A \}, \quad  a L_{G,\varphi,\psi}(A) \mapsto [a]
  \]
  is bijective, and hence
  $\lvert \{ [a] \in \mathcal{R}(\varphi,\psi) \mid a \in A \} \rvert
  = \lvert A : L_{G,\varphi,\psi}(A) \rvert$.
  % In particular, these considerations apply to $A = \overline{\gamma}_c(G)$.

  \smallskip

  (1) Suppose that $G$ has finite rank and that $R(\varphi,\psi) <
  \infty$.  By induction, the coincidence sets
  $\mathrm{Coin}(\varphi_i,\psi_i)$ are singletons for $1 \le i \le
  c-1$, and consequently 
  \[
  L_{G,\varphi,\psi}(A) = L_{A,\varphi_c,\psi_c}(A) \qquad \text{for
    $A = \overline{\gamma}_c(G)$.}
  \]
  For a contradiction, assume that $\mathrm{Coin}(\varphi_c,\psi_c)$
  is not a singleton.  Then the kernel of the endomorphism $A \to A$,
  $b \mapsto (b\varphi)^{-1} (b \psi)$ is non-trivial, hence its
  image, viz.\ $L_{A,\varphi_c,\psi_c}(A)$, has strictly smaller rank
  than~$A$.  (After tensoring the $\mathbb{Z}$-module $A$ with
  $\mathbb{Q}$, we are reduced to considering a linear endomorphism of
  the finite dimensional $\mathbb{Q}$-vector space
  $\mathbb{Q} \otimes_\mathbb{Z} A$.)  This implies that
  \[
  R(\varphi,\psi) \ge \lvert \{ [a] \in \mathcal{R}(\varphi,\psi) \mid
  a \in A \} \rvert = \lvert A : L_{A,\varphi_c,\psi_c}(A) \rvert = \infty.
  \] 

  \smallskip

  (2) Suppose that, for $1 \le i \le c$, the coincidence set
  $\mathrm{Coin}(\varphi_i,\psi_i)$ is a singleton.  We need to
  investigate how the pre-image of a twisted conjugacy class in
  $G / \overline{\gamma}_c(G)$ splits into twisted conjugacy classes
  in~$G$.  Fix $x \in G$ and consider
  $z, \tilde z \in \overline{\gamma}_c(G) \subseteq \mathrm{Z}(G)$.
  The elements $xz$ and $x\tilde z$ are $(\varphi,\psi)$-twisted
  conjugate if and only if there exists $g \in G$ such that
  ${\tilde z}^{-1} z = (g \vartheta)^{-1} (g \psi)$, where
  $\vartheta = \varphi \gamma_x$ and $\gamma_x \colon G \to G$,
  $h \mapsto x^{-1}hx$ denotes the inner automorphism that corresponds
  to conjugation by~$x$.

  Again we write $A = \overline{\gamma}_c(G)$, and we note that, for
  $1 \le i \le c$, the endomorphism
  $\vartheta_i \colon \overline{\gamma}_i(G) /
  \overline{\gamma}_{i+1}(G) \to \overline{\gamma}_i(G) /
  \overline{\gamma}_{i+1}(G)$
  induced by $\vartheta$ is equal to~$\varphi_i$.  In particular, the
  coincidence sets $\mathrm{Coin}(\vartheta_i,\psi_i)$, for
  $1 \le i \le c-1$, are singletons and
  $ L_{G,\vartheta,\psi}(A) = L_{A,\vartheta_c,\psi_c}(A)$.  Thus we
  conclude: the collection of twisted conjugacy classes in~$G$ which
  make up the pre-image of the twisted conjugacy class of the image of
  $x$ in $G / \overline{\gamma}_c(G)$ is in bijection to 
  \[
    A / L_{G,\vartheta,\psi}(A) = A / L_{A,\vartheta_c,\psi_c}(A) = A
    / L_{A,\varphi_c,\psi_c}(A) =
    \mathcal{R}(\varphi_c,\psi_c). \qedhere
  \]
 \end{proof}

\begin{Ex}
  It is easy to produce a torsion-free abelian group $A$, not of
  finite rank, and an endomorphism $\varphi \colon A \to A$ such that
  $R(\varphi) < \infty$, but with infinite coincidence set
  $\mathrm{Coin}(\varphi,\mathrm{id}_A) = \mathrm{Fix}(\varphi)$.  For
  instance, consider the
  Cartesian product $A = \prod_{i \in \mathbb{N}} \mathbb{Z}$ and the
  shift-map $\varphi \colon A \to A$,
  $(a_i)_{i \in \mathbb{N}} \mapsto (a_{i+1})_{i \in \mathbb{N}}$.  In
  this case, $R(\varphi)=1$, but
  $\mathrm{Coin}(\varphi,\mathrm{id}_A)$ consists of all constant
  sequences and is thus infinite.
\end{Ex}

\subsection{Profinite completions} \label{subsec:prof-compl}
Let $\varphi, \psi \colon G \to G$ be endomorphisms of a group~$G$.
This induces continuous endomorphisms
$\widehat{\varphi}, \widehat{\psi} \colon \widehat{G} \to \widehat{G}$
of the profinite completion~$\iota \colon G \to \widehat{G}$, and a
natural map
\[
\mathcal{R}(\varphi,\psi) \to
\mathcal{R}(\widehat{\varphi},\widehat{\psi}), \quad
[x]_{\varphi,\psi} \mapsto [x \iota
]_{\widehat{\varphi},\widehat{\psi}}.
\]
This constellation was already considered in \cite[Sec.~5.3.2]{Fe00}
and~\cite{FeTrZi20}.
 
\begin{Lem} \label{lem:prof-compl}
  In the situation described above, the following hold.
  \begin{enumerate}
  \item[$(1)$] Each $(\widehat{\varphi},\widehat{\psi})$-twisted
    conjugacy class in $\widehat{G}$ is compact and hence closed
    in~$\widehat{G}$.
  \item[$(2)$] For $x \in G$, the twisted conjugacy class
    $[x \iota ]_{\widehat{\varphi},\widehat{\psi}}$ is the closure of
    $[x]_{\varphi,\psi} \iota$ in~$\widehat{G}$.
  \item[$(3)$] If $R(\varphi,\psi) < \infty$, then the natural map
    $\mathcal{R}(\varphi,\psi) \to
    \mathcal{R}(\widehat{\varphi},\widehat{\psi})$ is surjective and
    each $(\widehat{\varphi},\widehat{\psi})$-twisted
    conjugacy class is open in~$\widehat{G}$.
  \item[$(4)$] If $G$ is abelian and $R(\varphi,\psi) < \infty$ then
    the natural map
    $\mathcal{R}(\varphi,\psi) \to
    \mathcal{R}(\widehat{\varphi},\widehat{\psi})$ is bijective.
  \end{enumerate}
\end{Lem}

\begin{proof}
  (1) For $x \in \widehat{G}$, the image of $\widehat{G}$ under the
  continuous map $\widehat{G} \to \widehat{G}$,
  $g \mapsto (g \widehat{\varphi})^{-1} x (g \widehat{\psi})$ is
  compact.

  (2) From (1) we conclude that
  $\overline{[x]_{\varphi,\psi} \iota} \subseteq [x \iota
  ]_{\widehat{\varphi},\widehat{\psi}}$.
  The reverse inclusion holds, because
  $\big( (g_\lambda \varphi)^{-1} x (g_\lambda \psi) \big) \iota$
  converges to
  $(g \widehat{\varphi})^{-1} (x \iota) (g \widehat{\psi})$,
  for any $g \in \widehat{G}$ and any net
  $(g_\lambda)_{\lambda \in \Lambda}$ in $G$ such that
  $(g_\lambda \iota)_{\lambda \in \Lambda}$ converges to $g$
  in~$\widehat{G}$.

  (3) Suppose that $R(\varphi,\psi) < \infty$.  Then every $y \in
  \widehat{G}$ is the limit of a net $(x_\lambda \iota)_{\lambda \in
    \Lambda}$, where the elements $x_\lambda \in G$ all belong to the
  same $(\varphi,\psi)$-twisted conjugacy class in~$G$.  The claim
  follows from~(2).

  (4) Suppose that $G$ is abelian and $R(\varphi,\psi) < \infty$.
  Then $\mathcal{R}(\varphi,\psi) = H_{\varphi,\psi} \backslash G$,
  where
  $H_{\varphi,\psi} = \{ (g\varphi)^{-1} (g\psi) \mid g \in G \} \le
  G$
  is a finite-index subgroup.  Also $\widehat{G}$ is abelian and
  $\mathcal{R}(\widehat{\varphi},\widehat{\psi}) =
  H_{\widehat{\varphi},\widehat{\psi}} \backslash G$,
  where
  $H_{\widehat{\varphi},\widehat{\psi}} = \{ (g\widehat{\varphi})^{-1}
  (g \widehat{\psi}) \mid g \in \widehat{G} \} \le \widehat{G}$
  equals the closure of $H_{\varphi,\psi} \iota$ in $\widehat{G}$ and
  $\lvert \widehat{G} : H_{\widehat{\varphi},\widehat{\psi}} \rvert =
  \lvert G : H_{\varphi,\psi} \rvert$.  The claim follows from~(3).
\end{proof}

\begin{Ex}
  For $G = \mathbb{Z}$ and $\varphi = \mathrm{id}_G$, the natural map
  $\mathcal{R}(\varphi) \to \mathcal{R}(\widehat{\varphi})$ is not
  surjective.  For $G = \mathbb{Q}$ and $\varphi = \mathrm{id}_G$, the
  natural map
  $\mathcal{R}(\varphi) \to \mathcal{R}(\widehat{\varphi})$ is not
  injective.
\end{Ex}

The following examples illustrate that an argument given in \cite[\S
5.3.2]{Fe00}, intended to show that every (finitely generated)
almost-abelian group contains a fully invariant, abelian finite-index
subgroup, is not quite complete.

\begin{Ex}
  Consider
  $G = \langle z \rangle \times (\langle a \rangle \ltimes \langle b
  \rangle ) \cong C_2 \times \mathrm{Sym}(3)$,
  where $\mathrm{ord}(z) = \mathrm{ord}(a) = 2$, $\mathrm{ord}(b) = 3$
  and $b^a = b^{-1}$.  The $2$-Sylow subgroup
  $A = \langle z,a \rangle \cong C_2 \times C_2$ is abelian and of
  index $\lvert G : A \rvert = 3$.  There are two other $2$-Sylow
  subgroups: $A^b = \langle z, ab^2 \rangle$ and
  $A^{b^2} = \langle z, ab \rangle$.  Thus the intersection of all
  subgroups of index $3$ in~$G$, namely $Z = \langle z \rangle$, is
  not fully invariant in~$G$; for instance, the endomorphism
  $\varphi \colon G \to G$ given by $z \mapsto a$ and $a,b \mapsto 1$
  does not map $Z$ to itself.
\end{Ex}

\begin{Ex}
  Consider
  $G = \langle b,a_0 \rangle \times \bigoplus_{i \in \mathbb{N}}
  \langle a_i \rangle \cong \mathrm{Sym}(3) \times \bigoplus_{i \in
    \mathbb{N}} C_2$, where $\mathrm{ord}(b) = 3$,
  $\mathrm{ord}(a_0) = \mathrm{ord}(a_1) = \ldots = 2$ and
  $b^{a_0} = b^{-1}$ , equipped with the endomorphism
  $\varphi \colon G \to G$ given by $b \varphi = a_0 \varphi = 1$ and
  $a_i \varphi = a_{i-1}$ for $i \in \mathbb{N}$.  Every
  $\varphi$-invariant finite-index subgroup of~$G$ contains an element
  of the form $c = a_i^{\, -1} a_j$ for $i,j \in \mathbb{N}$ with
  $i<j$ and hence $c \varphi^j = a_0$.  Consequently, the almost
  abelian group $G$ has no $\varphi$-invariant, abelian finite-index
  normal subgroup.
\end{Ex}

The following lemma partly fixes this small gap, in a slightly more
general context.

\begin{Lem}
  Let $G$ be finitely generated and almost-$\mathcal{P}$, where
  $\mathcal{P}$ is a group-theoretic property that is inherited by
  finite-index subgroups.  Then $G$ admits a fully invariant
  finite-index subgroup satisfying~$\mathcal{P}$.
\end{Lem}

\begin{proof}
  Put $d = \mathrm{d}(G)$, the minimal number of generators of~$G$,
  and let $A \le G$ be a subgroup that has property $\mathcal{P}$ and
  finite index $\lvert G : A \rvert < \infty$.  Replacing $A$ by its
  core in~$G$, we may assume that $A \trianglelefteq G$.  Denote by
  $W$ the set of group words $w = w(x_1,\ldots,x_d)$, i.\,e.\ elements
  in the free group on free generators $x_1, \ldots, x_d$, such that
  $w$ is a law in the finite group~$G/A$.  By~\cite[Thm.~15.4]{Ne67},
  the $d$-generated relatively free group in the variety corresponding
  to all laws in $G/A$ is finite.  Thus the verbal subgroup
  $B = \langle w(g_1,\ldots,g_d) \mid w \in W, \; g_1, \ldots, g_d \in
  G \rangle$ is fully invariant in $G$ and has finite index
  $\lvert G : B \rvert < \infty$.  As $B \le A$, we conclude that $B$
  satisfies~$\mathcal{P}$.
\end{proof}

The following proposition generalises Lemma~\ref{lem:prof-compl}(4)
and altogether by-passes some of the apparent difficulties in \cite[\S
5.3.2]{Fe00}.

\begin{Prop} \label{profinite-bijection}
  Let $\varphi, \psi \colon G \to G$ be endomorphisms of an almost
  abelian group~$G$ such that $R(\varphi,\psi) < \infty$.  
%   that is invariant under $\varphi$ and~$\psi$.
  Then the natural map
  $\mathcal{R}(\varphi,\psi) \to
  \mathcal{R}(\widehat{\varphi},\widehat{\psi})$ is a bijection.
\end{Prop}

\begin{proof}
  Fix $x \in G$.  In view of Lemma~\ref{lem:prof-compl}, it suffices
  to show that the map
  $\mathcal{R}(\varphi,\psi) \to
  \mathcal{R}(\widehat{\varphi},\widehat{\psi})$
  is injective.  For this it suffices to prove that, for every
  $x \in G$, there exists a finite-index subgroup $H \le G$ such that
  $[x]_{\varphi,\psi} \supseteq xH$; indeed, if $y \in G$ is such that
  $\overline{[x]_{\varphi,\psi} \iota} = \overline{[y]_{\varphi,\psi}
    \iota}$
  then $y$ is $(\varphi,\psi)$-twisted conjugate to an element of the
  open neighbourhood $x H$ of $x \in G$ (in the profinite topology
  on~$G$, which is not necessarily Hausdorff), hence $y$ is
  $(\varphi,\psi)$-twisted conjugate to~$x$.
 
  Fix $x \in G$, and let $A \le G$ be an abelian finite-index
  subgroup.  We observe that, for $a, \tilde a \in A$, the elements
  $xa$ and $x \tilde a$ are $(\varphi,\psi)$-twisted conjugate if and
  only if there exists $g \in G$ such that
  $a = (g \vartheta)^{-1} {\tilde a} (g \psi)$, where
  $\vartheta = \varphi \gamma_x$ and $\gamma_x \colon G \to G$,
  $h \mapsto x^{-1}hx$ denotes the inner automorphism that corresponds
  to conjugation by~$x$.  Put $A_0 = A \vartheta^{-1} \cap A \psi^{-1}
  \cap A$, a
  finite-index subgroup of~$A$, and let $Y$ be a right-transversal for
  $A_0$ in $G$ with~$1 \in Y$.  Writing $g = by$, with $b \in A_0$ and
  $y \in Y$, we obtain
  \begin{equation} \label{equ:umformung} (g \vartheta)^{-1} {\tilde a}
    (g \psi) = (y \vartheta)^{-1} \cdot {\tilde a} (b \vartheta)^{-1}
    (b \psi) \cdot (y \psi).
  \end{equation}
  We observe that
  $B = \{ (b \vartheta)^{-1} (b \psi) \mid b \in A_0 \}$ is a subgroup
  of~$A$.  Furthermore, \eqref{equ:umformung} shows that $x \tilde a$
  is $(\varphi,\psi)$-twisted conjugate to one of finitely many
  elements $x a_1, \ldots, x a_m$, where $a_1, \ldots, a_m \in A$, if
  and only if the coset $\tilde a B$ contains one of finitely many
  elements $(y \vartheta) a_i (y \psi)^{-1}$, with $y \in Y$ and
  $1 \le i \le m$.  Thus $R(\varphi,\psi) < \infty$ implies that
  $\lvert G : B \rvert = \lvert G : A \rvert \lvert A : B \rvert <
  \infty$.
  Taking $y=1$ in \eqref{equ:umformung}, we conclude that
  $[x]_{\varphi,\psi} \supseteq [x]_{\varphi,\psi} \cap xA \supseteq
  xB$.
\end{proof}

%%%%%

\section{An explicit formula for abelian groups of finite
  rank} \label{sec:explicit-formula}

Using techniques from commutative algebra, Miles derived
in~\cite{Mi08} a closed formula for periodic point counts for ergodic
finite-entropy endomorphisms of finite-dimensional compact abelian
groups.  This approach was then used in~\cite{BeMiWa14} to conjecture
a P\'olya--Carlson dichotomy for the analytic behaviour of dynamical
zeta functions of compact group automorphisms; supportive evidence
comes from a certain class of automorphisms of solenoids (connected
finite-dimensional compact abelian groups).  Under the Pontryagin
duality, finite-dimensional compact abelian groups correspond to
abelian groups of finite rank and, in the more restricted setting,
solenoids correspond to torsion-free abelian groups of finite rank.

Our aim in this section is to consider somewhat more general
situations and to arrive at a formula similar to Miles' by an entirely
different route, namely via profinite completions.  We start with a
reduction to torsion-free groups.

\begin{Lem}
  Let $G$ be an abelian group of finite rank, and let
  $\varphi, \psi \colon G \to G$ be endomorphisms such that
  $R(\varphi,\psi) < \infty$.  Let $T \le G$ denote the torsion
  subgroup of~$G$, and write
  $\tilde\varphi, \tilde\psi \colon G/T \to G/T$ and
  $\varphi_0, \psi_0 \colon T \to T$ for the induced endomorphisms.

  Then
  $R(\varphi,\psi) = R(\tilde\varphi,\tilde\psi) R(\varphi_0,\psi_0)$
  and, in particular,
  $R(\tilde\varphi,\tilde\psi), R(\varphi_0,\psi_0) < \infty$.
\end{Lem}

\begin{proof}
  We follow the same ideas as in the proof of Proposition~\ref{Nil}.
  Two elements $x,y \in G$ can be $(\varphi,\psi)$-twisted conjugate
  only if their images in $G/T $ are
  $(\tilde\varphi,\tilde\psi)$-twisted conjugate to one another.  Thus
  $R(\tilde\varphi,\tilde\psi) < \infty$, and
  $\mathrm{Coin}(\tilde\varphi,\tilde\psi)$ is a singleton by
  Proposition~\ref{Nil}.  We need to investigate how the pre-image of
  a twisted conjugacy class in $G/T$ splits into twisted conjugacy
  classes in~$G$.  Fix $x \in G$ and consider $z, \tilde z \in T$.
  The elements $xz$ and $x\tilde z$ are $(\varphi,\psi)$-twisted
  conjugate if and only if there exists $g \in G$ such that
  ${\tilde z}^{-1} z = (g \varphi)^{-1} (g \psi)$, because in contrast
  to the proof of Proposition~\ref{Nil} we can ignore inner
  automorphisms altogether.  Since
  $\mathrm{Coin}(\tilde\varphi,\tilde\psi)$ is trivial, we conclude
  that $xz$ and $x\tilde z$ are $(\varphi,\psi)$-twisted conjugate if
  and only if there exists $g \in T$ such that
  ${\tilde z}^{-1} z = (g \varphi)^{-1} (g \psi)$, which is to say
  that $z$ and $\tilde z$ are $(\varphi_0,\psi_0)$-twisted equivalent.
\end{proof}

\begin{Cor} \label{cor:reduction}
  Let $G$ be an abelian group of finite rank, and let
  $\varphi, \psi \colon G \to G$ be endomorphisms such that
  $(\varphi,\psi)$ is tame.  Let $T \le G$ denote the torsion subgroup
  of~$G$, and write $\tilde\varphi, \tilde\psi \colon G/T \to G/T$ and
  $\varphi_0, \psi_0 \colon T \to T$ for the induced endomorphisms.

  Then $(\tilde\varphi,\tilde\psi)$ and $(\varphi_0,\psi_0)$ are tame
  and
  \[
    Z_{\varphi,\psi}(s) = Z_{\tilde\varphi,\tilde\psi}(s) \ast
    Z_{\varphi_0,\psi_0}(s)
  \]
  is an additive convolution, in the sense of~\cite[\S 1]{FeHi94}.
\end{Cor}

\begin{Rmk}
  A similar result for finitely generated abelian groups and
  additional facts were established in \cite[Thm.~2]{FeHi94}.  The
  point here is that we may effectively replace the data
  $(G,\varphi,\psi)$ by the similar data
  $(G/T \times T, \tilde\varphi \times \varphi_0, \tilde\psi \times
  \psi_0)$ without any change to the coincidence Reidemeister zeta
  function.  Furthermore, by Proposition~\ref{profinite-bijection},
  the group $T$ can be replaced by $T_1 = T/T_\mathrm{div}$, where
  $T_\mathrm{div}$ denotes the maximal divisible subgroup of~$T$,
  equipped with the induced endomorphisms $\varphi_1, \psi_1$ of
  $T_1$: the group $T_1$ is of the form
  $\bigoplus_{p \in \mathbb{P}} F_p$, for finite abelian $p$-groups of
  uniformly bounded rank, and the profinite completion of $T$ is
  $\widehat{T} \cong \widehat{T_1} \cong \prod_{p \in \mathbb{P}}
  F_p$.  In the setting of \cite{Mi08}, the group $F_p$ is trivial for
  almost all primes~$p$; thus $T_1$ is finite and its
  $(\varphi_1,\psi_1)$-coincidence Reidemeister zeta function is
  easily understood; compare~\cite{FeHi94}.
\end{Rmk}

We now focus on the torsion-free case and derive a formula similar to
the one of Miles~\cite{Mi08}, but perhaps more directly accessible by
group-theoretic means.  For this we introduce the following notion of
`weak commutativity': We say that two endomorphisms
$\alpha, \beta \colon V \to V$ of a finite-dimensional
$\mathbb{Q}$-vector space $V$ are (absolutely) \emph{simultaneously
  triangularisable}, if there is a finite algebraic extension $L$ of
$\mathbb{Q}$ such that the induced endomorphisms of the $L$-vector
space $L \otimes_\mathbb{Q} V$ simultaneously preserve a complete flag
of subspaces.  For this to happen, it suffices that $\alpha$ and
$\beta$ commute; for a more precise characterisation of the property
we refer to~\cite{DrDuGr51}. 

\begin{Prop}\label{formula}
  Let $G$ be a torsion-free abelian group of finite rank~$d\geq 1$,
  and let $(\varphi,\psi)$ be a tame pair of endomorphisms
  for~$G$.  Let
  $\varphi_\mathbb{Q}, \psi_\mathbb{Q} \colon G_\mathbb{Q} \to
  G_\mathbb{Q}$ denote the extensions of $\varphi, \psi$ to the
  divisible hull
  $G_\mathbb{Q} = \mathbb{Q} \otimes_\mathbb{Z} G \cong \mathbb{Q}^d$
  of~$G$, and suppose that the endomorphisms
  $\varphi_\mathbb{Q}, \psi_\mathbb{Q}$ are simultaneously
  triangularisable.  Let $\xi_1, \ldots, \xi_d$ and
  $\eta_1, \ldots, \eta_d$ denote the eigenvalues of
  $\varphi_\mathbb{Q}$ and $\psi_\mathbb{Q}$ in a fixed algebraic
  closure of the field~$\mathbb{Q}$, including multiplicities,
  ordered so that, for $n \in \mathbb{N}$, the eigenvalues of
  $\varphi_\mathbb{Q}^{\, n} - \psi_\mathbb{Q}^{\, n}$ are
  $\xi_1^{\, n} - \eta_1^{\, n}, \ldots, \xi_d^{\, n} - \eta_d^{\,
    n}$.  Set
  $L = \mathbb{Q}(\xi_1, \ldots, \xi_d, \eta_1, \ldots, \eta_d)$.
 
  For each $p \in \mathbb{P}$, we fix an embedding
  $L \hookrightarrow \overline{\mathbb{Q}}_p$, where
  $\overline{\mathbb{Q}}_p$ denotes the algebraic closure of
  $\mathbb{Q}_p$, and we write $\lvert \cdot \rvert_p$ for the unique
  prolongation of the $p$-adic absolute value
  to~$\overline{\mathbb{Q}}_p$; this is a way of choosing a particular
  extension of the $p$-adic absolute value to~$L$.  Likewise, we
  choose an embedding $L \hookrightarrow \mathbb{C}$ and denote by
  $\lvert \cdot \rvert_\infty$ the usual absolute value
  on~$\mathbb{C}$.

  \smallskip

  Then there exist subsets $I(p) \subseteq \{1,\ldots,d\}$, for
  $p \in \mathbb{P}$, such that the following hold.
  \begin{enumerate}
  \item[$(1)$] For each $p \in \mathbb{P}$, the polynomials
    $\prod_{i \in I(p)} (X - \xi_i)$ and
    $\prod_{i \in I(p)} (X - \eta_i)$ have coefficients in
    $\mathbb{Z}_p$; in particular,
    $\lvert \xi_i \rvert_p , \lvert \eta_i \rvert_p \le 1$ for
    $i \in I(p)$.
  \item[$(2)$] For each $n \in \mathbb{N}$,
  \begin{equation}\label{equ:key-formula}
    R(\varphi^n,\psi^n) = \prod_{p \in \mathbb{P}} \prod_{i \in I(p)}
    \lvert \xi_i^{\, n} - \eta_i^{\, n} \rvert_p^{\, -1} =
    \prod_{i=1}^d \lvert \xi_i^{\, n} - \eta_i^{\, n} \rvert_\infty \cdot
    \prod_{p \in \mathbb{P}} \prod_{i \not \in I(p)} \lvert \xi_i^{\, n} -
    \eta_i^{\, n} \rvert_p; 
    \end{equation}
    as this number is a positive integer,
    $\lvert \xi_i^{\, n} - \eta_i^{\, n} \rvert_p = 1$ for
    $1 \le i \le d$ for almost all~$p \in \mathbb{P}$.
  \end{enumerate}
\end{Prop}

\begin{proof}
  Without loss of generality, $\mathbb{Z}^d \le G \le \mathbb{Q}^d$
  and we use additive notation.  Using Lemma~\ref{lem:prof-compl}(4), we
  may pass to the profinite completion.

  The profinite completions of the abelian group $G$ and its
  endomorphisms decompose as direct products
 \[
 \iota \colon G \to \widehat{G} = \prod_{p \in \mathbb{P}}
 \widehat{G}_p \qquad \text{and} \qquad \widehat{\varphi} = \prod_{p
   \in \mathbb{P}} \widehat{\varphi}_p, \quad \widehat{\psi} = \prod_{p
   \in \mathbb{P}} \widehat{\psi}_p,
 \]
 where, for each prime $p$, the Sylow pro-$p$ subgroup of
 $\widehat{G}$ constitutes the pro-$p$ completion $\widehat{G}_p$
 of~$G$, equipped with endomorphisms
 $\widehat{\varphi}_p, \widehat{\psi}_p \colon \widehat{G}_p \to
 \widehat{G}_p$.  Proposition~\ref{profinite-bijection} shows that
 \[
 R(\varphi^n,\psi^n) = \prod_{p \in \mathbb{P}}
 R(\widehat{\varphi}_p^{\, n}, \widehat{\psi}_p^{\, n}), \qquad
 \text{for $n \in \mathbb{N}$;}
 \]
 in particular, $R(\varphi^n,\psi^n) < \infty$ implies that
 $R(\widehat{\varphi}_p^{\, n}, \widehat{\psi}_p^{\, n}) = 1$ for
 almost all $p \in \mathbb{P}$ so that the product is only formally
 infinite.

 Fix $p \in \mathbb{P}$.  The pro-$p$ group $\widehat{G}_p$ is
 torsion-free, abelian and of rank at most~$d$, hence
 $\widehat{G}_p \cong \mathbb{Z}_p^{\, d(p)}$, where
 $d(p) = \mathrm{d}(\widehat{G}_p) \le d$.  Indeed, we can arrive at
 $\widehat{G}_p$ in the following way, which is somewhat roundabout
 but useful for our purposes.  We observe that the group $G$ has the
 same pro-$p$ completion as
 $G_{\mathbb{Z}_p} = \mathbb{Z}_p \otimes_\mathbb{Z} G$, where
 $\mathbb{Z}_p$ denotes the ring of $p$-adic integers.  The
 $\mathbb{Z}_p$-module $G_{\mathbb{Z}_p}$ decomposes as a direct sum
 of its maximal divisible submodule
 $D \cong \mathbb{Q}_p^{\, d - d(p)}$ and a free
 $\mathbb{Z}_p$-module~$A \cong \mathbb{Z}_p^{\, d(p)}$;
 compare~\cite[Thm.~20]{Ka54}.  Furthermore, the submodule $D$ is
 invariant under the induced $\mathbb{Z}_p$-module endomorphisms
 $\varphi_{\mathbb{Z}_p}, \psi_{\mathbb{Z}_p} \colon G_{\mathbb{Z}_p}
 \to G_{\mathbb{Z}_p}$;
 thus we obtain endomorphisms
 $\overline{\varphi}_{\mathbb{Z}_p}, \overline{\psi}_{\mathbb{Z}_p}
 \colon G_{\mathbb{Z}_p}/D \to G_{\mathbb{Z}_p}/D$.
 Consequently, there is an isomorphism
 $\sigma \colon \widehat{G}_p \to G_{\mathbb{Z}_p}/D$ of pro-$p$
 groups that is compatible with the endomorphism pairs
 $\widehat{\varphi}_p, \widehat{\psi}_p$ and
 $\overline{\varphi}_{\mathbb{Z}_p}, \overline{\psi}_{\mathbb{Z}_p}$,
 i.\,e.\
 \[
 \widehat{\varphi}_p \sigma = \sigma
 \overline{\varphi}_{\mathbb{Z}_p} \qquad \text{and}
 \qquad \widehat{\psi}_p \sigma = \sigma
 \overline{\psi}_{\mathbb{Z}_p}.
\]

We conclude that there exists a subset $I(p) \subseteq \{1,\ldots,d\}$
such that the endomorphisms $\widehat{\varphi}_p$ and
$\widehat{\psi}_p$ have eigenvalues $\xi_i$, $i \in I(p)$, and
$\eta_i$, $i \in I(p)$.  In particular, the polynomials
$\prod_{i \in I(p)} (X - \xi_i)$ and $\prod_{i \in I(p)} (X - \eta_i)$
have coefficients in $\mathbb{Z}_p$.

Finally, Lemma~\ref{lem:img-psi-central} yields
\[
  R(\widehat{\varphi}_p^{\, n},\widehat{\psi}_p^{\, n}) = \lvert
  \det(- \widehat{\varphi}_p^{\, n} + \widehat{\psi}_p^{\, n})
  \rvert_p^{\, -1} = \prod_{i \in I(p)} \lvert \xi_i^{\, n} -
  \eta_i^{\, n} \rvert_p^{\, -1}.
\]

Taking the product over all primes $p$, we arrive at the first formula
in~(2).  Using the adelic formula
$\lvert a \rvert_\infty \prod_{p \in \mathbb{P}} \lvert a \rvert_p =
1$
for $a \in \mathbb{Q} \smallsetminus \{0\}$, we obtain the second
formula in~(2).
\end{proof}

%%%%%

\section{P\'olya--Carlson dichotomy for coincidence Reidemeister zeta
  functions} \label{sec:polya-carlson}

In this section we present, in analogy to~\cite{BeMiWa14}, results in
support of a P{\'o}lya--Carlson dichotomy between rationality and a
natural boundary for the analytic behaviour of the coincidence
Reidemeister zeta functions for tame pairs of commuting endomorphisms
of a torsion-free nilpotent group of finite rank.  The classical
P{\'o}lya--Carlson theorem, as discussed in~\cite[\S 6.5]{Se08},
provides the following connection between the arithmetic properties of
the coefficients of a complex power series and its analytic behaviour.

\begin{PCThm}
  A power series with integer coefficients and radius of
  convergence~$1$ is either a rational function or has the unit circle
  as a natural boundary.
\end{PCThm}

Translated to our set-up, Bell, Miles and Ward~\cite{BeMiWa14}
conjectured that if $\varphi$ is a tame automorphism of an abelian
group $G$ of finite rank then $Z_\varphi(s)$ is either rational or
admits a natural boundary at its radius of convergence.  They
collected substantial evidence for this conjecture from a certain
class of automorphisms of torsion-free abelian groups of finite rank.

Let $G$ be a group with tame endomorphism pair $(\varphi,\psi)$.  It
is convenient to  introduce a notation
\[
  Z^*_{\varphi,\psi}(s) =\sum_{n= 1}^{\infty} R(\varphi^n,\psi^n)s^n =
  s \cdot Z'_{\varphi,\psi}(s)/Z_{\varphi,\psi}(s)
\]
for the generating series that enumerates directly the numbers of
coincidence Reidemeister classses.  The following lemma is basic: it
shows in particular that, if $Z^*_{\varphi,\psi}(s)$ has a natural
boundary at its radius of convergence, then so
does~$Z_{\varphi,\psi}(s)$; compare~\cite{BeMiWa14}.

\begin{Lem}\label{generating}
  If~$Z_{\varphi,\psi}(s)$ is rational then~$Z^*_{\varphi,\psi}(s)$ is
  rational.  If~$Z_{\varphi,\psi}(s)$ admits analytic continuation
  beyond its radius of convergence, then so
  does~$Z^*_{\varphi,\psi}(s)$.
\end{Lem}

For the proofs of the main theorems in this section we rely on the
following key result of Bell, Miles and Ward~\cite[Lem.~17]{BeMiWa14};
one of the ingredients in its proof is the Hadamard quotient theorem.

\begin{Lem}\label{derivative}
  Let~$S$ be a finite list of places of algebraic number fields and,
  for each~$v\in S$, let~$\xi_v$ be an element of the appropriate
  number field that is not a unit root and such
  that~$\lvert \xi_v \rvert_v=1$. Then the complex function
  \[
    F(s)=\sum_{n= 1}^{\infty}f(n)s^n, \qquad
    \text{where~$f(n)=\prod_{v\in S}\lvert \xi_v^n-1 \rvert_v$ for~$n\ge1$,}
  \]
  has the unit circle as a natural boundary.
\end{Lem}

First we state and prove our result in the case of abelian groups, in
order to illustrate some of the underlying ideas.

\begin{Thm} \label{abel} Let $\varphi,\psi \colon G \to G$ be a tame
  pair of endomorphisms of a torsion-free abelian group $G$ of finite
  rank~$d\geq 1$.  Using the notation from Proposition~\ref{formula},
  suppose that $\varphi_\mathbb{Q}$ and $\psi_\mathbb{Q}$ are
  simultaneously triangularisable and, in addition, that the primes
  $p$ that contribute non-trivial factors
  $\prod_{i \not\in I(p)} \lvert \xi_i^{\, n} - \eta_i^{\, n}
  \rvert_p$ to the product on the far right-hand side of
  \eqref{equ:key-formula} form a finite subset
  $\mathbf{P} \subseteq \mathbb{P}$.  For instance, this is the case,
  when $I(p) = \{1, \ldots, d\}$ for almost all $p \in \mathbb{P}$,
  equivalently when $G$ has no elements of infinite $p$-height for
  almost all $p \in \mathbb{P}$.  Suppose further that
  $\lvert \xi_i \rvert_\infty \neq \lvert \eta_i \rvert_\infty$ for
  $1 \le i \le d$.

  Then the coincidence Reidemeister zeta function $Z_{\varphi,\psi}(s)$
  is either a rational function or it has a natural boundary at its
  radius of convergence.  Furthermore, the latter occurs if and only if
  $\lvert \xi_i \rvert_p= \lvert \eta_i \rvert_p$ for
  some~$p \in \mathbf{P}$ and $i \not\in I(p)$.
\end{Thm}

\begin{proof}
  For $p \in \mathbf{P}$ we write
  \[
    S(p) = \{1,\ldots,d\} \smallsetminus I(p) \qquad \text{and} \qquad
    S^*(p) = \{i \in S(p) \mid \lvert \xi_i \rvert_p \neq \lvert
    \eta_i \rvert_p\};
  \]
  we remark right away that $\eta_i \neq 0$ for every
  $i \in S(p) \smallsetminus S^*(p)$, because otherwise
  $\xi_i = \eta_i = 0$ and $\varphi_\mathbb{Q} - \psi_\mathbb{Q}$
  would have rank less than $d$ contrary to
  Lemma~\ref{lem:img-psi-central}.  We set
  \[
    b = \prod_{p \in \mathbf{P}} \, \prod_{i \in S^*(p)} \max \{ \lvert
    \xi_i \rvert_p, \lvert \eta_i \rvert_p \}, \qquad
    \eta = \prod_{p \in \mathbf{P}} \, \prod_{i \in S(p)\smallsetminus
      S^*(p)}\lvert \eta_i \rvert_p
  \]
  and, for $n \in \mathbb{N}$,
  \[
    g(n) = b^n\cdot\eta^n \cdot \prod_{i=1}^d \lvert \xi_i^{\, n} - \eta_i^{\, n}
    \rvert_\infty  \qquad \text{and} \qquad
    f(n)=\prod_{p \in \mathbf{P}} \, \prod_{i \in S(p)\smallsetminus
      S^*(p)} \big\vert (\xi_i
      \eta_i^{\, -1})^n - 1 \big\vert_p.
  \]

  From formula~\eqref{equ:key-formula} in Proposition~\ref{formula} and using
  the ultrametric property, we deduce that, for every $n \in \mathbb{N}$,
  \begin{align*}
    R(\varphi^n,\psi^n) %
    & =  \prod_{i=1}^d \lvert \xi_i^{\, n} - \eta_i^{\, n} \rvert_\infty \cdot
      \prod_{p \in \mathbf{P}} \, \prod_{i  \in S(p)} \lvert \xi_i^{\, n} -
      \eta_i^{\, n} \rvert_p \\
    % & = \prod_{i=1}^d \lvert \xi_i^{\, n} - \eta_i^{\, n}
    %   \rvert_\infty \cdot \prod_{p \in \mathbf{P}} \prod_{i \in S^*(p)}
    %   \lvert \xi_i^{\, n} -  \eta_i^{\, n} \rvert_p\prod_{i \in
    %   S(p)\smallsetminus S^*(p)}\lvert \xi_i^{\, n} -  \eta_i^{\, n}
    %   \rvert_p \\
    & = \prod_{i=1}^d \lvert \xi_i^{\, n} - \eta_i^{\, n}
      \rvert_\infty \cdot \prod_{p \in \mathbf{P}} \, \prod_{i \in S^*(p)}
      \max\{\lvert \xi_i \rvert_p^n, \lvert \eta_i \rvert_p^n\} \cdot
      \prod_{p \in \mathbf{P}} \, \prod_{i \in S(p)\smallsetminus
      S^*(p)}\lvert \xi_i^{\, n} - \eta_i^{\, n} \rvert_p \\   
    & = \prod_{i=1}^d \lvert \xi_i^{\, n} - \eta_i^{\, n}
      \rvert_\infty \cdot b^n \cdot \eta^n \cdot \prod_{p \in
      \mathbf{P}} \, \prod_{i \in S(p)\smallsetminus S^*(p)} \big\vert (\xi_i
      \eta_i^{\, -1})^n - 1 \big\vert_p \\
    & = g(n) f(n).
  \end{align*}

  Now we open up the absolute values in the product
  $P(n) = \prod_{i=1}^d \lvert \xi_i^{\, n} - \eta_i^{\, n}
  \rvert_\infty$.  Complex eigenvalues $\xi_i$ in the spectrum of
  $\varphi_\mathbb{Q}$, respectively $\eta_i$ in the spectrum of
  $\psi_\mathbb{Q}$, appear in pairs with their complex conjugate
  $\overline{\xi_i}$, respectively $\overline{\eta_i}$.

  Moreover, such pairs can be lined up with one another in a
  simultaneous triangularisation as follows. Write $\varphi_L, \psi_L$
  for the induced endomorphisms of the $L$-vector space
  $V = L \otimes_\mathbb{Q} G \cong L^d \hookrightarrow \mathbb{C}^d$.
  If $v \in V$ is, at the same time, an eigenvector of $\varphi_L$
  with complex eigenvalue $\xi_d$ and an eigenvector of $\psi_L$ with
  eigenvalue $\eta_d$, then there is $w \in V$ such that $w$ is, at
  the same time, an eigenvector of $\varphi_L$ with eigenvalue
  $\overline{\xi_d} \ne \xi_d$ and an eigenvector of $\psi_L$ with
  eigenvalue $\overline{\eta_d}$, possibly equal to~$\eta_d$.  Thus we
  can start our complete flag of $\{\varphi,\psi\}$-invariant
  subspaces of $V$ with
  $\{0\} \subset \langle v \rangle \subset \langle v,w \rangle$, and
  proceed with $V/\langle v,w \rangle$ by induction to produce the
  rest of the flag in the same way, treating complex eigenvalues of
  $\psi_L$ in the same way as they appear.

  If at least one of $\xi_i, \eta_i$ is complex so that these
  eigenvalues of $\varphi_\mathbb{Q}$ and $\psi_\mathbb{Q}$ are paired
  with eigenvalues
  $\xi_j = \overline{\xi_i}, \eta_j = \overline{\eta_i}$, for suitable
  $j \ne i$, as discussed above, we see that
    \[
      \big\lvert \xi_i^{\, n} - \eta_i^{\, n} \big\rvert_\infty \;
      \big\lvert \xi_j^{\, n} - \eta_j^{\, n} \big\rvert_\infty =
      \big\lvert \xi_i^{\, n} - \eta_i^{\, n} \big\rvert_\infty ^{\,
        2} = \big( \xi_i^{\, n} - \eta_i^{\, n} \big) \cdot
      (\overline{\xi_i}^{\, n} - \overline{\eta_i}^{\, n} \big).
    \]
    If $\xi_i$ and $ \eta_i$ are both real eigenvalues of
    $\varphi_\mathbb{Q}$ and $\psi_\mathbb{Q}$, not paired up with
    another pair of eigenvalues, then exactly as in
    Example~\ref{ex1} above we have
    $\lvert \xi_i^{\, n} - \eta_i^{\, n} \rvert_\infty=
    \delta_{1,i}^{\, n} - \delta_{2,i}^{\, n}$, where
    $\delta_{1,i} = \max\{\lvert \xi_i \rvert_\infty,\lvert \eta_i
    \rvert_\infty\} $ and
    $\delta_{2,i}=\frac{\xi_i\cdot\eta_i}{\delta_{1,i}}$.

  Hence we can expand the product $P(n)$, using an appropriate
  symmetric polynomial, to obtain an expression of the form
  \begin{equation}\label{dominant1}
    g(n) = \sum_{j \in J} c_jw_j^{\, n},
  \end{equation}
  where~$J$ is a finite index set, $c_j \in \{-1,1\}$ and
  $\{ w_j \mid j \in J \} \subseteq \mathbb{C} \smallsetminus \{0\}$.
  Consequently, the coincidence Reidemeister zeta function can be
  written as
  \[
    Z_{\varphi,\psi}(s) = \exp \left( \sum_{j \in J} c_j \sum_{n=
        1}^{\infty} \frac{f(n) (w_js)^n} {n} \right).
  \]
  If~$S(p) = S^*(p)$ for all~$p\in \mathbf{P}$, then~$f(n) = 1$ for
  all $n \in \mathbb{N}$ and it follows that
  $Z_{\varphi,\psi}(s) = \prod_{j \in J} (1 - w_js)^{-c_j}$ is a
  rational function.

  Now suppose that $S(p) \neq S^*(p)$ for some~$p\in \mathbf{P}$.  By
  Lemma~\ref{generating}, it suffices to show that 
  \[
    Z^*_{\varphi,\psi}(s) = \sum_{j \in J} c_j\sum_{n= 1}^\infty f(n)
    (w_js)^n
  \]
  has a natural boundary at its radius of convergence.  Moreover,
  from~
  Lemma~\ref{derivative} we have
  $\limsup_{n\rightarrow\infty}f(n)^{1/n}=1$.  Hence, for
  each~$j \in J$, the series
  \[
    \sum_{n= 1}^{\infty} f(n)(w_js)^n
  \]
  has radius of convergence~$\lvert w_j\lvert_\infty^{-1}$.

  As $\lvert \xi_i \rvert_\infty \neq \lvert \eta_i \rvert_\infty$ for
  $1 \le i \le d$, there is a dominant term~$w_m$ in the
  expansion~\eqref{dominant1} for which
  \[
    \lvert w_m \lvert_\infty = b\cdot\eta \cdot \prod_{i=1}^d
    \max\{\lvert \xi_i \rvert_\infty,\lvert \eta_i \rvert_\infty\}
  \]
  and~$\lvert w_m \lvert_\infty > \lvert w_j \lvert_\infty$ for
  all~$j \in J \smallsetminus \{m\}$.  Thus it suffices to show that
  $\sum_{n= 1}^{\infty} f(n)(w_m s)^n$ has its circle of convergence
  as a natural boundary.  This is the case, because
  $\sum_{n= 1}^{\infty} f(n)s^n$ has the unit circle as a natural
  boundary by Lemma~\ref{derivative}.
\end{proof}

We remark that the special case of Theorem~\ref{abel}, where
$\psi = \mathrm{id}_G$ and $\varphi$ is an automorphism of the abelian
group $G$ has been discussed in \cite[Thm.~5.8]{FeZi20}, via methods
from commutative algebra as in~\cite{BeMiWa14}.  The following example
of a coincidence Reidemeister zeta function with natural boundary
comes from a $3$-adic extension of the circle doubling map and is
essentially known; see~\cite{EvStWa05}.

\begin{Ex}
  Consider the endomorphisms $\varphi \colon g \mapsto 6g$ and
  $\psi \colon g \mapsto 3g$ of the abelian group $\Z[\frac{1}{3}]$,
  written additively.  A straightforward calculation, in line with
  Lemma~\ref{formula}, shows that, for $n \in \mathbb{N}$,
  \[
    R(\varphi^n, \psi^n)= \lvert \coker(\varphi^n-\psi^n) \lvert  =
    \lvert 6^n-3^n \rvert_{\infty}\cdot \lvert 6^n-3^n \rvert_3 = 
    \lvert 2^n-1 \rvert_{\infty}\cdot \lvert 2^n-1 \rvert_3. 
  \]
  Hence, the coincidence Reidemeister zeta function is
  \[
    Z_{\varphi,\psi}(s) = \exp \left( \sum_{n=1}^\infty \frac{\lvert
        2^n-1 \rvert_{\infty}\cdot \lvert 2^n-1 \rvert_3 }{n} s^n
    \right),
  \]
  This zeta function was studied in~\cite{EvStWa05}.  Their calculations
  show that the modulus of the coincidence Reidemeister zeta function
  satisfies
  \begin{equation*}
    \lvert Z_{\varphi,\psi}(s) \rvert_\infty =
    \left\vert\frac{1-s}{1-2s}\right\vert_\infty \cdot
    \left\vert\frac{1-(2s)^{2\vphantom{3^j}}}{1-s^2}\right\vert_\infty^{\frac{1}{2}}
    \cdot
    \left\vert\frac{1-s^{2\vphantom{3^j}}}{1-(2s)^{2}}\right\vert_\infty^{\frac{1}{6}}
    \cdot
    \prod_{j=1}^{\infty}\left\vert
      \frac{1-(2s)^{2\cdot 3^j}}{1-s^{2 \cdot 3^j}}\right\vert_\infty^{1/(2
      \cdot 9^j)}.
  \end{equation*}
  It follows that the series defining the zeta function has zeros at
  all points of the form $\frac{1}{2}e^{2\pi i /3^r}$, $r \ge1$,
  whence $\vert s\vert=\frac{1}{2}$ is a natural boundary for the
  coincidence Reidemeister zeta function $Z_{\varphi,\psi}(s)$.
\end{Ex}

Our proof of Theorem~\ref{abel} requires, at a technical level, the
finiteness of the set $\mathbf{P} \subseteq \mathbb{P}$ and the
assumption that
$\lvert \xi_i \rvert_\infty \neq \lvert \eta_i \rvert_\infty$ for
$1 \le i \le d$, but there is no indication that these conditions are
actually necessary for the desired dichotomy to hold; regarding the
second condition, compare the comment following Theorem~15
in~\cite{BeMiWa14}.  In the next example we look at a coincidence
Reidemeister zeta function for a tame pair of endomorphisms of
$G \cong \mathbb{Z}^2$ that cannot be simultaneously triangularised;
the example shows that a possible P\'olya--Carlson dichotomy in this
generality needs to allow for new outcomes.

\begin{Ex} \label{exa:non-comm-endos}
  Consider the automorphisms
  $\varphi, \psi \colon \mathbb{Z}^2 \to \mathbb{Z}^2$ that are given
  by
  \[
    (x,y) \varphi = (x,y) \begin{pmatrix} 1 & 1 \\ 0 & 1 \end{pmatrix}
    = (x,x+y), \quad (x,y) \psi = (x,y) \begin{pmatrix} 1 & 0 \\ 1 &
      1 \end{pmatrix} = (x+y,y).
  \]
  Using Lemma~\ref{lem:img-psi-central}, it is easy to check that, for
  $n \in \mathbb{N}$,
  \[
    R(\varphi^n, \psi^n)= \lvert \coker(\varphi^n-\psi^n) \lvert =
    n^2.
  \]
  Hence, the coincidence Reidemeister zeta function is
  \[
    Z_{\varphi,\psi}(s) = \exp \left( \sum_{n=1}^\infty n s^n
    \right) = \exp \left( \frac{s}{(1-s)^2}\right).
  \]
  Clearly, $Z_{\varphi,\psi}(s)$ is neither a rational function nor
  does the zeta function have a natural boundary at its radius of
  convergence.
\end{Ex}

Finally we generalise Theorem~\ref{abel} to deduce the main result
that was stated in the introduction.

\begin{proof}[Proof of Theorem~\ref{thm:main-result-nilpotent}]
  Assertion (1) follows from Proposition~\ref{Nil}, and (2) follows
  from Proposition~\ref{formula}.  Thus it remains to prove (3), and
  we follow closely the proof of Theorem~\ref{abel}.

  For $1 \le k \le c$ and $p \in \mathbf{P}$ we write
  $S_k(p) = \{1, \ldots, d_k\} \smallsetminus I_k(p)$ and
  $S^*_k(p) = \{i \in S_k(p) \mid \lvert \xi_{k,i} \rvert_p \neq
  \lvert \eta_{k,i} \rvert_p\}$.  We set
  \[
    b = \prod_{k=1}^c \prod_{p \in \mathbf{P}_k} \prod_{i \in
      S^*_k(p)} \max \{ \lvert \xi_{k,i} \rvert_p, \lvert \eta_{k,i} \rvert_p
    \}, \qquad \eta = \prod_{k=1}^c \prod_{p \in \mathbf{P}_k} \, \prod_{i \in
      S_k(p) \smallsetminus S^*_k(p)}\lvert \eta_{k,i} \rvert_p,
  \]
  and as in the proof of Theorem~\ref{abel} we deduce from (1) and (2)
  that, for $n \in \mathbb{N}$,
  \[
    R(\varphi^{\, n},\psi^{\, n}) = g(n) \cdot f(n),
  \]
  where
  \[
    g(n)= b^n \cdot \eta^n \cdot \prod_{k=1}^c \prod_{i=1}^{d_k}
    \lvert \xi_{k,i}^{\, n} - \eta_{k,i}^{\, n} \rvert_\infty \quad
    \text{and} \quad f(n)= \prod_{k=1}^c \prod_{p \in \mathbf{P}_k}
    \prod_{i \in S_k(p) \smallsetminus S^*_k(p)} \big\vert (\xi_{k,i}
    \eta_{k,i}^{\, -1})^n - 1 \big\vert_p.
  \]
  As in the proof of Theorem \ref{abel} we open up the absolute values
  in the product
  $P(n) = \prod_{k=1}^c \prod_{i=1}^{d_k} \lvert \xi_{k,i}^{\, n} -
  \eta_{k,i}^{\, n} \rvert_\infty$ to arrive at an expression of the
  form
  \[
    g(n) = \sum_{j \in J} c_jw_j^{\, n},
  \]
  where~$J$ is a finite index set, $c_j \in \{-1,1\}$ and
  $\{ w_j \mid j \in J \} \subseteq \mathbb{C} \smallsetminus \{0\}$.
  Consequently, the coincidence Reidemeister zeta function can be
  written as
  \[
    Z_{\varphi,\psi}(s) = \exp \left( \sum_{j \in J} c_j \sum_{n=
        1}^{\infty} \frac{f(n) (w_js)^n} {n} \right).
  \]
  and we conclude the argument along the same line as in the proof of
  Theorem~\ref{abel}.
\end{proof}

% \noindent
% \textbf{Acknowledgement.} We thank the referee for valuable feedback
% which led to improvements in the exposition.

%%%%% 
%%%%%


\begin{thebibliography}{10}

\bibitem{BeMiWa14} J.~Bell, R.~Miles, and T.~Ward, \textit{Towards a
    P\'olya-Carlson dichotomy for algebraic dynamics}, Indag.\ Math.\
  (N.S.) \textbf{25} (2014), 652--668.

\bibitem{ByCo18} J.~Byszewski and G.~Cornelissen, \textit{Dynamics on
    abelian varieties in positive characteristic}, with an appendix by
  R.~Royals and T.~Ward,  Algebra Number Theory \textbf{12} (2018), 
  2185--2235.

\bibitem{DeDu15} K.~Dekimpe and G.-J.~Dugardein, \textit{Nielsen zeta
    functions for maps on infra-nilmanifolds are rational}, J.\ Fixed
  Point Theory Appl.\ \textbf{17} (2015), 355--370.

\bibitem{DeTeBu18} K.~Dekimpe, S.~Tertooy, and I.~Van~den Bussche,
  \textit{Reidemeister zeta functions of low-dimensional
    almost-crystallographic groups are rational}, Comm.\ Algebra
  \textbf{46} (2018), 4090--4103.

\bibitem{DrDuGr51} M.\,P.\ Drazin, J.\,W.\ Dungey, and K.\,W.\
  Gruenberg, \textit{Some theorems on commutative matrices}, J.\
  London Math.\ Soc.\ \textbf{26} (1951), 221--228.

\bibitem{EvStWa05} G.~Everest, V.~Stangoe, and T.~Ward, \textit{Orbit
    counting with an isometric direction}, 293--302, in:
  \emph{Algebraic and topological dynamics}, Contemp.\ Math.\
  \textbf{385}, Amer.\ Math.\ Soc., Providence, RI, 2005.
  
\bibitem{Fe91} A.~L.~Fel'shtyn, \textit{The Reidemeister zeta function and the
  computation of the Nielsen zeta function}, Colloq.\ Math.\ \textbf{62} (1991),
  153--166.

\bibitem{Fe00} A.~Fel'shtyn, \textit{Dynamical zeta functions, Nielsen
    theory and Reidemeister torsion}, Mem.\ Amer.\ Math.\ Soc.\
  \textbf{147} (2000), no. 699.

% \bibitem{Fe10} A.~Fel'shtyn, \textit{New directions in
%     Nielsen-Reidemeister theory}, Topology Appl.\ 157 (2010),
%   no. 10-11, 1724--1735.

% \bibitem{FeGo08} A.~Fel'shtyn, D.\,L.~Gon\c{c}alves, \textit{The
%     Reidemeister number of any automorphism of a Baumslag--Solitar
%     group is infinite}, in: Geometry and dynamics of groups and
%   spaces, 399--414, Progr.\ Math., 265, Birkh\"auser, Basel, 2008.

\bibitem{FeHi94} A.~Fel'shtyn and R.~Hill, \textit{The Reidemeister
    zeta function with applications to Nielsen theory and a connection
    with Reidemeister torsion}, $K$-Theory \textbf{8} (1994),
  367--393.

\bibitem{FeLe15} A.~Fel'shtyn and  J.\,B.~ Lee, \textit{The Nielsen and
     Reidemeister numbers of maps on infra-solvmanifolds of type
    $(R)$}, Topology Appl.\ \textbf{181} (2015), 62--103.

% \bibitem{FeLeTr08} A.~Fel'shtyn, Y.~Leonov, E.~Troitsky,
%   \textit{Twisted conjugacy classes in saturated weakly branch
%     groups}, Geom.\ Dedicata 134 (2008), 61--73.

% \bibitem{FeNa16} A.~Fel'shtyn, T.~Nasybullov, \textit{The $R_\infty$
%     and $S_\infty$ properties for linear algebraic groups}, J.\ Group
%   Theory 19 (2016), no. 5, 901--921.

% \bibitem{FeTr08} A.~Fel'shtyn, E.~Troitsky, \textit{Geometry of
%     Reidemeister classes and twisted Burnside theorem}, J. K-Theory 2
%   (2008), no. 3, 463--506.

% \bibitem{FeTr18} A.~Fel'shtyn, E.~Troitsky, \textit{Twisted
%     Burnside-Frobenius theory for endomorphisms of polycyclic groups},
%   Russ.\ J.\ Math.\ Phys.\ 25 (2018), no. 1, 17--26.

% \bibitem{FeTrVe06} A.~Fel'shtyn, E.~Troitsky, A.~Vershik,
%   \textit{Twisted Burnside theorem for type$II_1$ groups: an example},
%   Math.\ Res.\ Lett.\ 13 (2006), no. 5-6, 719--728.

\bibitem{FeTrZi20} A.~Fel'shtyn, E.~Troitsky, and M.~Zietek,
  \textit{New zeta functions of Reidemeister type and twisted
    Burnside-Frobenius theory}, Russ.\ J.\ Math.\ Phys.\ \textbf{27}
  (2020),199--211.

\bibitem{FeZi20} A.~Fel'shtyn and M.~Zietek, \textit{Dynamical zeta
    functions of Reidemeister type and representations spaces},
  57--81, in: Contemp.\ Math.\ \textbf{744}, Amer.\ Math.\ Soc.,
  Providence, RI, 2020.

\bibitem{Ka54} I.~Kaplansky, Infinite abelian groups, University of
  Michigan Press, Ann Arbor, 1954.

% \bibitem{Le07} G.~Levitt, \textit{On the automorphism group of
%     generalized Baumslag-Solitar groups}, Geom.\ Topol.\ 11 (2007),
%   473--515.

\bibitem{Mi08} R.~Miles, \textit{Periodic points of endomorphisms on
    solenoids and related groups}, Bull.\ Lond.\ Math.\ Soc.\
  \textbf{40} (2008), 696--704.
 
\bibitem{Mi13} R.~Miles, \textit{Synchronization points and
    associated dynamical invariants}, Trans.\ Amer.\ Math.\ Soc.\
  \textbf{365} (2013), 5503--5524.

% \bibitem{Na12} T.\,R.~Nasybullov, \textit{Twisted conjugacy classes in
%     general and special linear groups} Algebra Logic 51 (2012), no. 3,
%   220--231.

\bibitem{Ne67} H.~Neumann, Varieties of groups, Springer-Verlag, New
  York, 1967.

% \bibitem{RiZa10} L.~Ribes, P.~Zalesskii, Profinite groups,
%   Springer-Verlag, Berlin, 2010.
  
\bibitem{Pa77} D.~Passman, The algebraic structure of group
  rings, Wiley Interscience (John Wiley and Sons), New York, 1977.

\bibitem{Ro96} D.\,J.\,S.\ Robinson, A course in the theory of groups,
  Springer-Verlag, New York, 1996.

\bibitem{Ro11} V.~Roman'kov, \textit{Twisted conjugacy classes in
    nilpotent groups}, J.\ Pure Appl.\ Algebra \textbf{215} (2011),
  664--671.
  
\bibitem{Se08} S.~L. Segal, Nine introductions in complex analysis,
  Elsevier Science B.V., Amsterdam, revised edition, 2008.

% \bibitem{Ru62} W.~Rudin, Fourier analysis on groups, Interscience
%   Publishers (a division of John Wiley and Sons), New York 1962.

% \bibitem{TrXX} E.~Troitsky, \textit{Reidemeister classes in some
%     weakly branch groups}, preprint (2018), \texttt{arXiv:1804.03695}.

\bibitem{We73} B.\,A.\,F.\ Wehrfritz, Infinite linear groups,
  Springer-Verlag, New York, 1973.

\bibitem{Wo01} P. Wong, \textit{Reidemeister zeta function for group
  extensions}, J.\ Korean Math.\ Soc.\ \textbf{38} (2001), 1107--1116.
\end{thebibliography}
\end{document}